\documentstyle[11pt]{article}

\def\b#1{{\bf #1}}
\def\i#1{{\it #1}}

\def\e#1{{\rm e}^{#1}}
\def\m#1{\mit #1}

\begin{document}

\title{Sampling from a couple of negatively correlated gamma variates}
\author{Mario Catalani
\\Department of Economics, University of Torino\\
Via Po 53, 10124 Torino, Italy\\
E-mail mario.catalani@unito.it}
\date{}
\maketitle
\begin{abstract}
\small{We propose two algorithms for sampling from two gamma variates
possessing a negative correlation. The case of positive correlation is
easily solved, so we just mention it. The main problem is the lowest value
of the correlation coefficient that can be reached. The starting point of
both algorithms
is generation from a bivariate density with uniform negatively correlated
marginals. Actually the first method uses a degenerate bivariate density
since it considers two uniforms related by a linear relationship.
Then we resort
essentially to the inverse transform method. For both algorithms we stress
restrictions on the parameters and rigidities.
}
\end{abstract}

\section{Introduction}
To fix notation, if $X$ follows a gamma distribution with parameters
$\lambda$ and $\alpha$, that is $X\sim {\cal G}(\lambda,\,\alpha)$,
then the density is
$$f_X(x)={\lambda^{\alpha}\e{-\lambda x}x^{\alpha-1}\over
{\m{\Gamma}}(\alpha)},\qquad x>0.$$
Johnson and Kotz \cite[1972]{kotz}
define a multivariate gamma distribution
in the following way. We have $m+1$ independent standard (that is,
each with $\lambda =1$), gamma variates (in general with different $\alpha's$)
$X_0,\,X_1,\,\ldots,\,X_m$. Define
$$Y_j=X_0+X_j, \qquad j=1,\,\ldots,\,m$$
Then the $Y_j,\,j=1,\,\ldots,\,m$ are distributed according a $m$-variate
gamma variable. We can see that each marginal $Y_j$ is a standard
gamma variable with $\alpha = \alpha_0 +\alpha_j$ and
\begin{eqnarray*}
\b{C}ov(Y_j,\,Y_{j'})&=& \b{V}ar(X_0)\\
&=&\alpha_0,
\end{eqnarray*}
and so all the marginals are positively correlated.

\noindent
We follow this suggestion, with $m=2$, but we exploit the preserving
monotonicity property of
the inverse transform method, see Fishman \cite[1996]{fishman}
to generate a couple of gamma
variates, $X_1$ and $X_2$, with a rigid amount of negative correlation
and then generate $X_0$ to allow for some flexibility.

\section{First Method}
Let $X_1\sim {\cal G}(1,\,r)$ and $X_2\sim {\cal G}(1,\,s)$, with $r$ and
$s$ integers. Assume, without loss of generality, $r<s$.
Let $U_i$ be $s$ independent uniforms on the unit interval, that is
$U_i\sim {\cal U}(0,\,1),\,i=1,\,\ldots,\,s$. Then we can write, using
properties of the gamma variate and the inverse transform method,
$$X_1=-\sum_{i=1}^r\ln U_i,$$
$$X_2=-\sum_{i=1}^s\ln (1-U_i).$$
It follows that
\begin{eqnarray*}
\b{C}ov(X_1,\,X_2)&=&\b{C}ov\left (\sum_{i=1}^r\ln U_i,\,
\sum_{i=1}^s\ln (1-U_i)\right )\\
&=&\sum_{i=1}^r\b{C}ov\left (\ln U_i,\, \ln (1-U_i)\right )\\
&=&\sum_{i=1}^r\left [\left (2-{\pi^2\over 6}\right )-1\right ]\\
&=&r\left (1-{\pi^2\over 6}\right ).
\end{eqnarray*}
Now define, with $X_0\sim {\cal G}(1,\,\alpha_0)$ independent of $X_1$
and $X_2$,
$$Y_1=X_0+X_1,$$
$$Y_2=X_0+X_2.$$
It follows $Y_1\sim {\cal G}(1,\,\alpha_0+r)$ and
$Y_2\sim {\cal G}(1,\,\alpha_0+s)$. Furthermore
\begin{eqnarray*}
\b{C}ov(Y_1,\,Y_2)&=&\b{V}ar(X_0)+\b{C}ov(X_1,\,X_2)\\
&=&\alpha_0 + r\left (1-{\pi^2\over 6}\right ).
\end{eqnarray*}
Because
$$c= 1-{\pi^2\over 6}=-0.644934,$$
we can see that varying $r$ and $\alpha_0$ we have some freedom in modelling
a negative covariance.

\noindent
The correlation coefficient between $Y_1$ and $Y_2$ is given by
\begin{equation}
\label{eq:correl}
\rho(Y_1,\,Y_2)={\alpha_0+rc\over\sqrt{(\alpha_0+r)(\alpha_0+s)}}.
\end{equation}
For $\alpha_0+rc<0$, which guarantees a negative correlation, we see
that
$${\partial \rho(Y_1,\,Y_2)\over \partial \alpha_0}>0,$$
and so when the correlation is negative, it is an increasing function of
$\alpha_0$. The lower bound for $\rho$ is given by
$${rc\over \sqrt{rs}}=c\sqrt{{r\over s}}.$$
When $r=s$ and assuming that $\alpha_0=0$ means that $X_0=0$ with
probability 1 we reach the most negative correlation possible between
two gamma variates.

\noindent
If we are given $\alpha_0+r=m$, $\alpha_0+s=n$ and $\rho=\rho_0$ then
from Equation~\ref{eq:correl} we get
\begin{equation}
\label{eq:correla}
\rho_0={m+r(c-1)\over\sqrt{mn}},
\end{equation}
which can be solved for $r$ and subsequently obtaining $\alpha_0$ and $s$.
The drawback is that $r$ and $s$ are required to be integers, so we
cannot arbitrarily choose $m$, $n$ and $\rho$. To give an idea of
the situation we present some tables, assuming for simplicity that
both $m$ and $n$ are integers. It follows that also $\alpha_0$ has
to be an integer. Because we want a negative correlation and because
$\alpha_0 \ge 0$ we have the restriction
$$r\le m < r(1-c).$$
Using Equation~\ref{eq:correla} we obtain

\bigskip

\begin{tabular}{||lllc||lllc||lllc||}\hline
{\em r}&{\em m}&{\em s}&{$\rho$}&
{\em r}&{\em m}&{\em s}&{$\rho$}&
{\em r}&{\em m}&{\em s}&{$\rho$}\\ \hline
2&2&3&-.5266&2&2&5&-.4078&2&2&8&-.3224\\
2&3&4&-.0837&2&3&6&-.0683&2&3&9&-.0557\\
5&5&6&-.5887&5&5&8&-.5098&5&5&11&-.4348\\
5&6&7&-.3432&5&6&9&-.3027&5&6&12&-.2621\\
5&7&8&-.1636&5&7&10&-.1463&5&7&13&-.1283\\
5&8&9&-.0264&5&8&11&-.0239&5&8&14&-.0212\\
8&8&9&-.6080&8&8&11&-.5500&8&8&14&-.4875\\
8&9&10&-.4384&8&9&12&-.4002&8&9&15&-.3579\\
8&10&11&-.3012&8&10&13&-.2771&8&10&16&-.2497\\
8&11&12&-.1879&8&11&14&-.1740&8&11&17&-.1579\\
8&12&13&-.0928&8&12&15&-.0864&8&12&18&-.0788\\
8&13&14&-.0118&8&13&16&-.0110&8&13&19&-.0101\\
12&12&13&-.6196&12&12&15&-.5768&12&12&18&-.5265\\
12&13&14&-.4995&12&13&16&-.4672&12&13&19&-.4288\\
12&14&15&-.3960&12&14&17&-.3720&12&14&20&-.3429\\
12&15&16&-.3059&12&15&18&-.2884&12&15&21&-.2670\\
12&16&17&-.2267&12&16&19&-.2144&12&16&22&-.1993\\
12&17&18&-.1565&12&17&20&-.1485&12&17&23&-.1385\\
12&18&19&-.0940&12&18&21&-.0894&12&18&24&-.0836\\
12&19&20&-.0379&12&19&22&-.0361&12&19&25&-.0338
\\ \hline
\end{tabular}

\bigskip

\noindent
Nothing essentially changes if we require gamma's with the first parameter
different from 1:
it is enough to multiply $Y_1$ and $Y_2$ by,
say, a constant $\beta$. The correlation coefficient is of course not
affected.

\section{Joint Density}
Obtaining the joint density of $Y_1$ and $Y_2$ in the general case is
cumbersome. To give an idea we consider the simplest case,
$r=s=1$, and we set $\alpha < -c$ that is $\alpha < 0.644934$, to guarantee
a negative correlation. So we have
\begin{equation}
Y_1=X_0-\ln U,
\end{equation}
\begin{equation}
Y_2=X_0-\ln (1-U).
\end{equation}
with independence of $X_0$ and $U$.
It turns out that marginally $Y_1$ and $Y_2$ are both
${\cal G}(1,\,1+\alpha_0)$.
Setting
$$\left\{\begin{array}{lll}y_1&=&x_0-\ln u\\
y_2&=&x_0-\ln (1-u),\end{array}\right . $$
we get
$$\left\{\begin{array}{lll}u&=&{1\over 1+\e{y_1-y_2}}\\
x_0&=&y_1-\ln \left (1+\e{y_1-y_2}\right ).\end{array}\right . $$
Because $X_0\ge 0$ we have the inequality
$$y_1-\ln \left (1+\e{y_1-y_2}\right )\ge 0.$$
The solution to the equation
$$y_1-\ln \left (1+\e{y_1-y_2}\right )= 0$$
is given by
$$y_2=y_1-\ln \left (\e{y_1}-1\right ).$$
It follows that we have $x_0 \ge 0$ if
$$y_2\ge y_1-\ln \left (\e{y_1}-1\right ).$$
For the Jacobian of this transformation we have
$$J={\e{y_1-y_2}\over \left (1+\e{y_1-y_2}\right )^2}.$$
The joint density of $X_0$ and $U$ is given by
$$f_{X_0,U}(x_0,\,u)=\e{-x_0}x_0^{\alpha -1}, \qquad \left \{\begin{array}{ll}
0\le & u\le 1\\
0\le & x_0. \end{array}\right . $$
The quantity $\e{-x_0}J$ simplifies to
$${1\over \e{y_1}+\e{y_2}}$$
and consequently we get for
the joint density of $Y_1$ and $Y_2$
$$f_{Y_1,Y_2}(y_1,\,y_2)=
{\left [y_1-\ln \left (1+\e{y_1-y_2}\right )\right ]^{\alpha -1}
\over \e{y_1}+\e{y_2}}, \qquad
\left \{\begin{array}{ll}
y_1 > & 0\\
y_2 > & y_1-\ln \left (\e{y_1}-1\right ). \end{array}\right . $$
We see that for $\alpha =1$ this simplifies to
$$f_{Y_1,Y_2}(y_1,\,y_2)=
{1\over \e{y_1}+\e{y_2}}, \qquad
\left \{\begin{array}{ll}
y_1 > & 0\\
y_2 > & y_1-\ln \left (\e{y_1}-1\right ), \end{array}\right . $$
but in this case the correlation is positive:
\begin{eqnarray*}
\rho(Y_1,\,Y_2)&=&{1+c\over 2}\\
&=&0.177533.
\end{eqnarray*}
In this last situation the marginals are  ${\cal G}(1,\,2)$. This can be
checked using the fact that
$$\int_b^{+\infty}{1\over a+\e{x}}\,dx ={-b+\ln\left (a+\e{b}\right )\over
a}.$$

\section{Second Method}
Again looking for flexibility we analyze another method where we start
generating $s$ samples $(U_{1_i},\,U_{2_i})$ from a bivariate distribution
with density
\begin{equation}
\label{eq:density2}
f(u_1,\,u_2)=1+\theta(1-2u_1)(1-2u_2),\qquad 0\le u_1,\, u_2\le 1.
\end{equation}
This density is of the form studied by Long and Krzysztofowicw
\cite[1995]{long}. For
$-1\le \theta\le 1$ this function is a proper density. It turns out
that marginally $U_{1_i}$ and $U_{2_i}$ are uniforms over the unit
interval and the correlation coefficient is given by
$$\rho(U_{1_i},\,U_{2_i})={\theta\over 3}.$$
Now define
$$X_1=-\sum_{i=1}^r\ln U_{1_i},$$
$$X_2=-\sum_{i=1}^s\ln U_{2_i}.$$
Now
\begin{eqnarray*}
\b{C}ov(\ln U_{1_i},\,\ln U_{2_i})&=& \int_0^1\int_0^1\ln x\ln y
(1+\theta (1-2x)(1-2y))\,dxdy -1\\
&=&1+{\theta\over 4} -1\\
&=&{\theta\over 4}.
\end{eqnarray*}
It follows that
$$
\b{C}ov(X_1,\,X_2)={r\theta\over 4}.$$
Proceeding as in the other method
we define, with $X_0\sim {\cal G}(1,\,\alpha_0)$ independent of $X_1$
and $X_2$,
$$Y_1=X_0+X_1,$$
$$Y_2=X_0+X_2,$$
and we get
\begin{equation}
\rho(Y_1,\,Y_2)={\alpha_0+{r\theta\over 4}\over
\sqrt{(\alpha_0+r)(\alpha_0+s)}}.
\end{equation}
Now we have a negative correlation if $4\alpha_0+r\theta<0$. Because of this
inequality the admissible range for $\theta$ is $-1\le\theta <0$. The
lower bound for $\rho$ is
$${\theta\over 4}\sqrt{r\over s},$$
so we cannot hope to do better than
$$\rho>-0.25.$$
\noindent
Repeating the same reasoning that led us to Equation~\ref{eq:correla} now we
have
$$\rho_0={4m-r(4-\theta)\over 4\sqrt{mn}}.$$
The restrictions are now
\begin{equation}
\label{eq:restr}
r\le m<{r(4-\theta)\over 4}.
\end{equation}
Solving for $r$ we get
\begin{equation}
\label{eq:erre}
r={4m-4\rho_0\sqrt{mn}\over 4-\theta}.
\end{equation}
Other restrictions on the minimum value admissible for $\rho_0$ arise
from Equation~\ref{eq:restr} and the fact that $r$ must be integer.
The minimum value of $r$ as given by Equation~\ref{eq:erre} is obtained
for $\theta=-1$. Then
$${4m-4\rho_0\sqrt{mn}\over 5}\le m-1,$$
which implies
\begin{equation}
\label{eq:soluzione}
\rho_0\ge -{m-5\over 4\sqrt{mn}}.
\end{equation}
Under the condition $m\ge 6$ this attainable lower bound is an increasing
function of $n$ and because $n\ge m$ it is a decreasing function of $m$:
in particular if $m=n$ the limit for $m$ going to infinity is $-0.25$.

\noindent
Once observed the conditions on $m$ and $\rho$ we have the following
algorithm.
\begin{enumerate}
\item Set $y=4m-4\rho_0\sqrt{mn}$.
\item Obtain $a={y\over 5}$.
\item Set $r=\lceil a\rceil$, that is the lowest integer not lower than $a$.
Because of the construction such an $r$ exists and $r<m$.
\item Obtain $\theta = 4-{y/r}$.
\end{enumerate}
As an example, imagine $m=7,\,n=10$. The attainable lower bound is
$-0.0597$. Set for example $\rho_0=-0.05$. Then ${y\over 5}=5.2653$. Take
$r=6$, so that $\theta =4-4.38778=-0.38778$.

\noindent
Once $\theta$ is obtained, the next step is to evaluate $\alpha_0=m-r$
and $s=n-\alpha_0$. Then we generate $s$ samples from the density
given in Equation~\ref{eq:density2}: this can be done for instance with
the acceptance-rejection method. And then we follow the final construction
to obtain $Y_1$ and $Y_2$.

\noindent
We can note that in this second method the attainable lower bound
for $\rho$ is sensibly
greater than in the first method, but for the admissible values we have a
complete flexibility. The other negative point is that we have to generate
samples from a bivariate density, whose \i{covariance scaler} $\theta$
(as termed by Long and Krzysztofowicw) has to be preliminarily obtained,
instead of generations from univariate densities.

\end{document}